\documentclass[11pt,twoside]{article}

\usepackage{latexsym, amssymb, amsmath, theorem}
\usepackage{epsfig}

\numberwithin{equation}{section}

\theoremstyle{change}
\theorembodyfont{\itshape}
\newtheorem{theorem}{Theorem}[section]
\newtheorem{proposition}[theorem]{Proposition}
\newtheorem{lemma}[theorem]{Lemma}
\newtheorem{corollary}[theorem]{Corollary}

\theorembodyfont{\rmfamily}
\newtheorem{definition}[theorem]{Definition}

\newtheorem{example}[theorem]{Example}
\newtheorem{remark}[theorem]{Remark}
\newenvironment{proof}{{\noindent \textbf{Proof}\,\,}}{\hspace*{\fill}$\Box$\medskip}

\newcounter{amarker}
\renewcommand{\theamarker}{\textup{(\arabic{amarker})}}

\newcounter{rmarker}

\makeatletter
\renewcommand{\@makecaption}[1]{%
\begin{center}#1\end{center}%
}
\makeatother
 
\title{\textbf{On the monodromy group of confluenting linear equations}}
\author{Alexey GLUTSYUK\\[5pt]}

\begin{document}
\maketitle 

\centerline{{\it Dedicated to Yu.S.Ilyashenko on the occasion of his 60$^{\text th}$ 
birthday.}}

\def\Dwidetilde{\dot{\mathstrut{\widetilde{\mathstrut z}}}} 
\def\var{\varepsilon} 
\def\im{\operatorname{Im}} 

\begin{abstract}
We consider a linear analytic ordinary differential equation with complex time 
having a nonresonant irregular singular point. We study it as a limit of a 
generic family of equations with confluenting Fuchsian singularities.

In 1984 V.I.Arnold asked the following question: is it true that some operators  
from the  monodromy group 
of the perturbed (Fuchsian) equation tend to Stokes operators of the 
nonperturbed irregular equation? 
Another version of this question was also independently proposed by J.-P.Ramis 
in 1988. 

We consider the case of Poincar\'e rank 1 only. 
We show (in dimension two) that generically no monodromy 
operator tends to a Stokes operator; on the other hand, in any dimension 
commutators of appropriate noninteger powers of 
the monodromy operators around singular points tend to Stokes operators.
\end{abstract}

\def\la{\lambda}

2000 Math. Subj. Class. 34M35 (34M40)


\

\centerline{\Large{ Table of contents.}}

\

1. Introduction.

1.1. Brief statements of results, the plan of the paper 
and the history.

1.2. Analytic classification of irregular equations. 
Canonical solutions and Stokes operators. 

1.3. Historical overview.

\

2. Main results. Stokes operators and limit monodromy. 

2.1. Stokes operators as limit transition operators between 
monodromy eigenbases.

2.2. Stokes operators as limits of commutators of appropriate 
powers of the monodromy operators.

2.3. The case of higher Poincar\'e rank.  

\

3. Convergence of the commutators to Stokes operators. 
Proof of Theorem 2.16.

3.1. Properties of the monodromy and the transition operators. 
The plan of the proof of Theorem 2.16.

3.2. The upper-triangular element of the transition matrix. Proof 
of (3.3). 

3.3. Commutators of operators with asymptotically common eigenline. 
Proof of Lemma 3.5. 

3.4. Convergence of the commutators to Stokes operators: 
the higher-dimensional case.

\

4. Generic divergence of monodromy operators along degenerating loops. 
 
 4.1 The statement of the divergence Theorem.
 
 4.2. Projectivization. The scheme of the proof of Theorem 4.6.
 
 4.3. Divergence of word of projectivizations. Proof of Lemma 4.10. 
 
 \
 
 References.

\section{Introduction} 

\subsection{Brief statements of results, the plan of the paper  
and the history} 

Consider a linear analytic ordinary differential equation 
\begin{equation} 
\dot z=\frac{A(t)}{t^{k+1}}z,\quad z\in\mathbb C^n,\,\,|t|\leq1,\,\,k\in 
\mathbb N\tag*{(1.1)}
\end{equation} 
with a nonresonant irregular singularity of order (the Poincar\'e rank) $k$ 
at 0 (or briefly, an irregular equation).  This means that 
$A(t)$ is a holomorphic matrix  
function such that the matrix $A(0)$ has distinct eigenvalues 
(denote them by $\lambda_i$). Then the matrix $A(0)$ is diagonalizable, and 
without loss of generality we suppose that it is diagonal.   
\def\a{\alpha} 
\begin{definition}  Two equations of type (1.1) are analytically (formally) 
equivalent, if there exists a change $z=H(t)w$ of the variable 
$z$, where $H(t)$ is a holomorphic invertible 
matrix function (respectively, a formal invertible matrix power series),  
that transforms one equation into the other. 
\end{definition}

The analytic classification of irregular equations (1.1) is well-known 
([2], [3], [9], [10], [16]): the complete system of invariants for 
analytic classification consists of a formal normal form (1.4) and 
Stokes operators (1.6) defined in 1.2; the latters are linear 
operators acting in the solution space of (1.1) comparing appropriate 
"sectorial canonical solution bases". 

On the other hand, an irregular equation (1.1) can be regarded as a result 
of {\it confluence of Fuchsian singular points} 
(recall that a Fuchsian singular point of a linear equation is a first 
order pole of its right-hand side). Namely, consider a deformation 
\begin{equation}
\dot
z=\frac{A(t,\var)}{f(t,\var)}z,	\ \ \ f(t,\var)=\prod_{i=0}^{k}(t-\a_i(\var)),
\tag{1.2}
\end{equation}
of equation (1.1) 
that splits the irregular singular point 0 of the nonperturbed equation into 
$k+1$ Fuchsian singularities $\a_i(\var)$ of the perturbed equation, i.e., 
$\a_i(\var)\neq\a_j(\var)$ for $i\neq j$.  
The family (1.2) depends on a parameter $\var\in \mathbb R_+\cup0$, 
$f(t,0)\equiv t^{k+1}$, $A(t,0)\equiv A(t)$.  

The {\it monodromy group} of a Fuchsian equation 
acts linearly in its solution space by 
analytic extensions of solutions along closed loops.
The analytic equivalence class of a 
%corr
generic 
%corr
Fuchsian equation 
is completely determined by the local types of its singularities and 
the action of its monodromy group. 

Everywhere below by $M_i$ we denote the monodromy operator of the perturbed 
equation (1.2) along a loop going around the singular point 
$\a_i$ (the choice of the corresponding loops will be specified later). 
The monodromy group of the perturbed equation is generated by appropriately 
chosen operators $M_i$. 

In 1984 V.I.Arnold proposed the following 
question. Consider a generic deformation (1.2). Is there an operator 
\begin{equation}
M_{i_l}^{d_l}\dots M_{i_1}^{d_1}\tag{1.3}
\end{equation}
from the monodromy group of the perturbed equation that converges to a Stokes 
operator of the nonperturbed equation? 

A version of this question was proposed independently by J.-P.Ramis in 1988. 

 \def\jlp{\text{[2], [3], [9], [10], [16]}}
It appears that already in the simplest case of dimension 2 and Poincar\'e 
rank $k=1$ generically 
{\it each operator from the monodromy group 
(except for that along a circuit (and its powers) around both singularities) 
tends to infinity} (Theorem 4.6 in Section 4), so, no one tends to a Stokes operator. 

In other terms, generically, no word (1.3) 
with $d_i\in\mathbb Z$ tends to a Stokes operator. But if $k=1$, then 
appropriate words (1.3) with {\it noninteger} exponents $d_i$ tend to Stokes 
operators (Theorem 2.16 in 2.2). 

The previous question and its nonlinear analogues  
were studied by J.-P.Ramis, B.Khesin, A.Duval, C.Zhang,   
J.Martinet, the author and others (see the historical overview in 1.3).
It was proved by the author in [6]  
that appropriate branches of the eigenfunctions of the monodromy operators 
$M_i$ of the perturbed equation tend to appropriate canonical solutions of the 
nonperturbed equation (Theorem 2.5 in 2.1). In the case of Poincar\'e rank 
$k=1$ this implies (Corollary 2.6 in 2.1) 
that Stokes operators of the nonperturbed equation are 
limits of transition operators between appropriate eigenbases of the 
monodromy operators $M_i$. This Corollary has a generalization 
for higher Poincar\'e rank [6]. 
  
The proofs of the results of the present paper are based on the previously 
mentioned results from [6], which are recalled in 2.1.

 In 1.2 we recall the analytic classification of irregular equations (1.1) 
 and the definitions of sectorial canonical solution bases and Stokes operators.  

  In 2.2 we state Theorem 2.16 on convergence of appropriate word (1.3) with 
  noninteger exponents $d_i$ to a Stokes operator in the case of Poincar\'e 
  rank $k=1$. Its proof is given in Section 3. 
  The corresponding exponents $d_i$  do not depend on the choice of 
  deformation. In fact, in the case of the higher Poincar\'e rank $k=2$ and 
  $n=2$ one can also prove a similar statement, but now the choice of the 
  corresponding exponents $d_i$ will depend on the choice of deformation. The 
  latter case will be discussed in 2.3.  
 
  In Section 4 in the case, when $k=1$, $n=2$, 
  for a typical nonperturbed equation (1.1) we prove the divergence of the  
operators from the monodromy group of the perturbed equation 
(except for the monodromy along a circuit around both singularities and its 
powers).  

\subsection{Analytic classification of irregular equations. 
Canonical solutions and Stokes operators}

Let (1.1) be an irregular equation, $\lambda_i$, $i=1,\dots,n$, be the eigenvalues 
of the corresponding matrix $A(0)$. 
 
One can	ask the question: is it true that the
variables $z=(z_1,\dots,z_n)$	in the equation
can be separated, more precisely, that (1.1) is analytically equivalent to a
direct sum of one-dimensional linear equations,	i.e.,  a
linear equation	with a diagonal	matrix function
in the right-hand side?	 Generically, the answer is ``no''.  At	the
same time any irregular equation (1.1) is formally equivalent to a unique 
direct sum of the type 
\begin{equation}
\begin{cases} \dot{w_i}=\frac{b_i(t)}{t^{k+1}}\,w_i,\\
i=1,\dots,n,
\end{cases}\tag{1.4}
\end{equation}
where $b_i(t)$ are polynomials of degrees at most $k$,
$b_i(0)=\lambda_i$.  The normalizing series bringing (1.1) to (1.4) 
is unique up to left
multiplication by constant diagonal matrix.  The system (1.4) is
called the {\it formal normal form} of (1.1) ($\jlp$).
\medskip

\def\re{\operatorname{Re}}
\def\wt#1{\widetilde#1}

Generically the	normalizing series diverges.  At the same time there
exists a finite covering
$\bigcup\limits_{j=0}^NS_j$ of a punctured neighborhood
of zero	in the $t$- line 
by radial sectors $S_j$ (i.e., those with the vertex at 0) that have
the following property. There exists a unique change of variables
$z=H_j(t)w$ over each $S_j$ that transforms (1.1) to $(1.4)$, where $H_j(t)$ 
is an analytic invertible matrix	
function on $S_j$ that can be $C^\infty$-smoothly extended to the
closure	$\overline S_j$	of the sector so that its asymptotic Taylor
series at 0 coincides with the normalizing series. The previous statement on 
existence and uniqueness of sectorial normalization holds in any good sector 
(see the two following Definitions); the covering consists of good sectors 
($\jlp$).

{\bf Case $k=1$, $n=2$, $\lambda_1-\lambda_2\in\mathbb R$.}

\begin{definition} A sector in $\mathbb C$ with the vertex at 0 is said to be 
{\it good}, if it contains only one imaginary semiaxis $i \mathbb R_{\pm}$, and 
its closure does not contain the other one (see Fig.1).
\end{definition}

{\bf General case.}

\begin{definition} (see, e.g., [9]). Let $k\in\mathbb N$, 
$\Lambda=\{\lambda_1,\dots,\lambda_n\}\subset\mathbb C$ be a $n$- ple of 
distinct numbers,  $t$ be the 
coordinate on $\mathbb C$. For a given pair $\lambda_i\neq\lambda_j$ 
the rays in $\mathbb C$ starting at 0 and 
forming the set    $\re(\frac{\lambda_j-\lambda_i}{t^k})=0$
are called the $(k,\Lambda)$- {\it imaginary dividing rays} corresponding
to the pair $(\lambda_i,\lambda_j)$. A radial sector 
is said to be  {\it $(k,\Lambda)$- good}, if for any pair
$(\lambda_i,\lambda_j)$, $j\neq i$, it contains exactly one corresponding
imaginary dividing ray and so does its closure.
\end{definition}

\begin{remark} In the case, when $k=1$, $n=2$, 
$\lambda_1-\lambda_2\in\mathbb R$, 
the imaginary dividing rays are the imaginary semiaxes, and the notions of 
good sector and $(k,\Lambda)$- good sector coincide. 
\end{remark}

\begin{remark} The ratio $\frac{w_i}{w_j}(t)$ of solutions of equations
from (1.4) tends to either zero or infinity, as $t$ tends to zero along a
ray distinct from the imaginary dividing rays corresponding to the
pair $(\lambda_i,\lambda_j)$. Its limit changes exactly when the
ray under consideration jumps over one of the latter
imaginary dividing rays.
\end{remark}

\def\diag{\operatorname{diag}}
\medskip
We consider a covering $\bigcup\limits_{j=0}^NS_j$ of a punctured neighborhood 
of zero by good (or $(k,\Lambda)$- good) 
sectors numerated counterclockwise and put $S_{N+1}=S_0$.
The standard splitting of the normal form  (1.4) into the direct sum
of one-dimensional equations defines a canonical base in	its solutions
space (uniquely	up to multiplication of	the base functions by constants)
with a diagonal fundamental matrix. Denote the latter fundamental
matrix by
$$W(t)=\diag(w_1,\dots,w_n).$$
Together with the normalizing changes $H_j$ in $S_j$, it defines the
canonical bases $(f_{j1},\dots,f_{jn})$ in the solution space of (1.1) in
the sectors $S_j$ with the fundamental matrices
\begin{equation}
Z^j(t)=H_j(t)W(t), \ j=0,\dots,N+1,\tag{1.5}\end{equation}
\noindent 
where for any $j=0,\dots,N$ the
branch ("with the index $j+1$") of the fundamental matrix $W(t)$ in $S_{j+1}$
is obtained from that
in $S_j$ by the counterclockwise analytic extension for any $j=0,\dots,N$.
(We put $S_{N+1}=S_0$. The corresponding branch of $W$ "with the index
$N+1$" is obtained from that "with the index 0" by the right multiplication
by the monodromy matrix of the formal normal form (1.4).) In a  connected 
component of the intersection $S_j\cap S_{j+1}$ there are two canonical solution 
bases coming from $S_j$ and $S_{j+1}$. Generically, they do not coincide.  
The transition between them is defined by a constant matrix $C_j$: 
\begin{equation}Z^{j+1}(t)=Z^j(t)C_j.\tag{1.6}\end{equation}
 The \ transition\ operators (matrices $C_j$)
are called {\it Stokes operators (matrices)} (see [2], [3], [9], 
[10], [16]).	The
nontriviality of Stokes	operators yields the obstruction to analytic
equivalence  of (1.1) and its formal normal form (1.4).
\medskip
\noindent
\begin{remark} Stokes matrices (1.6) are well-defined up to simultaneous
conjugation by one and the same	diagonal matrix.
\end{remark}
\medskip
\noindent
\begin{example}  Let $k=1,n=2$. In this case without
loss of generality we
assume that $\lambda_1-\lambda_2\in \mathbb R_+$ (one can achieve this by
linear change of the time variable). Then
the above covering consists of two sectors $S_0$ and $S_1$
(Figure	1).  The former	contains the positive imaginary	semiaxis and
its closure does not contain the negative one; the latter has the
same properties	with respect to	the negative (respectively, positive)
imaginary semiaxis.
There are two components of the intersection $S_0\cap S_1$.  So, in
this case we have a pair of Stokes operators.  {\it The
Stokes matrices (1.6) are
unipotent:} the one corresponding to the left intersection component 
 is	lower-triangular; the other one
is upper-triangular ($\jlp$).
\end{example}

\begin{figure}[ht]
  \begin{center}
   \epsfig{file=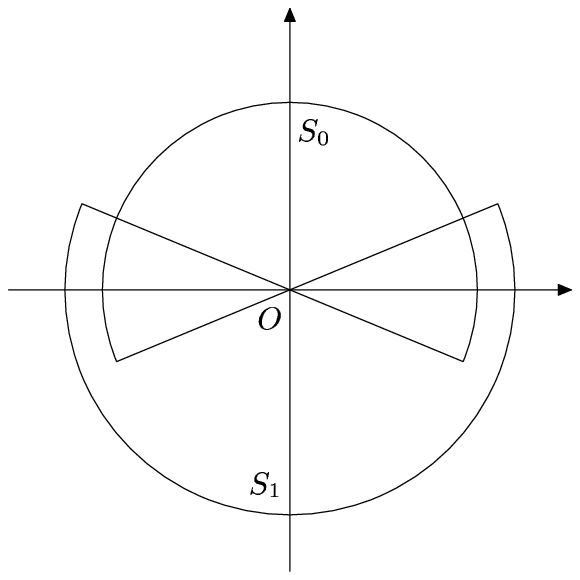}
    \caption{}
    \label{fig:1}
  \end{center}
\end{figure}

\begin{remark} Stokes operators of an irregular equation (1.1) with a diagonal
matrix in the right-hand side are identity operators. In this case (1.1) is
analytically equivalent to its formal normal form. In general, 
{\it two irregular equations are analytically
equivalent, if and only if they have the same formal normal form and
the corresponding Stokes matrix	tuples are obtained from each other by
simultaneous conjugation by one	and the	same diagonal matrix}, cf.
the previous Remark.  Thus, formal normal form and Stokes matrix tuple taken up
to the previous conjugation present the complete system of
invariants for analytic	classification of irregular equations (see $\jlp$).
\end{remark}

\subsection{Historical overview}

   Earlier in	 1919	R.   Garnier  [5]  had	studied	 some
particular deformations	of some	class of  linear equations with
nonresonant   irregular   singularity.  He obtained some 
analytic classification	invariants for these equations by studying their
     deformations.    The    complete	system	  of   analytic
classification	invariants   (Stokes  operators	and formal normal form)
for   general irregular  differential
equations was obtained later in	1970-th	years in  the papers by
Jurkat,	Lutz, Peyerimhoff [10], Sibuya [16] and Balser,
Jurkat,	 Lutz	[3]. Later Jurkat, Lutz and Peyerimhoff had extended 
their results to some resonant cases [11]. It is well-known that the monodromy 
operators of a linear ordinary differential equation belong to its Galois group 
(see [9], [14]). In 1985 J.-P.Ramis have proved that the Stokes operators also 
belong to the Galois group ([14], see also [9]). 
In 1989 he considered the classical confluenting
family of  hypergeometric  equations and proved	convergence of appropriate
branches of monodromy eigenfunctions of	the perturbed equation to
canonical solutions of the nonperturbed one by  direct calculation [15]. 
In the late 1980-th years B.Khesin also proved a
version	of  this  statement, but his result  was  not
published. In 1991  A.Duval  [4] proved this statement for the 
biconfluenting family  of  hypergeometric  equations (where the
nonperturbed equation  is  equivalent  to  Bessel  equation) by
direct calculation. In	1994 C.	 Zhang	[17] had	obtained  the
expression of  Garnier's  invariants  via  Stokes operators (for the 
class of irregular equations considered by Garnier). 

The conjecture saying that Stokes operators are limit transition operators 
between  monodromy eigenbases of the perturbed equation was firstly 
proposed by A.A.Bolibrukh in   1996. It was proved by the author in [6]. Later 
this result was extended to a generic resonant case [8]. 

Nonlinear\ analogues\ \ of\ the\ previous\ statements\ \ for\ parabolic\  
mappings \ 
(i.e., one-dimensional conformal mappings tangent to identity) and their 
\'Ecalle-Voronin moduli, saddle-node 
singularities of two-dimensional holomorphic vector fields and their 
Martinet-Ramis invariants (sectorial central manifolds in higher dimensions) 
 were obtained by the author in [7]. 
 Generalizations and other versions of the statement on parabolic mappings 
 were recently obtained in the joint paper [12] by P.Mardesic, R.Roussarie, 
 C.Rousseau, and in two unpublished joint papers 
 by the following authors: 1) X.Buff and Tan Lei; 2) A.Douady, 
 Francisco Estrada, P.Sentenac.

A particular case of the result from [7] concerning parabolic 
mappings (analogous to the previously mentioned statements 
on linear equations) was obtained by J.Martinet [13].

\section{Main results. Stokes operators and limit monodromy}
Everywhere below (whenever the contrary is not specified) we consider that 
the (nonperturbed) irregular equation under consideration has Poincar\'e 
rank $k=1$. 
In the present Section we recall the statements from [6] expressing 
the Stokes operators as limit transition operators between monodromy 
eigenbases of the confluenting Fuchsian equation (Theorems 2.5, 2.11 and 
Corollary 2.6 in Subsection 2.1). In 2.2 we state the results expressing the 
Stokes operators 
as limits of some words (1.3) of noninteger powers of monodromy 
operators (Theorem 2.16). In 2.3 we discuss the extension 
of these results to the case of higher Poincar\'e rank 
%corr
in dimension two. 
%corr

\subsection{Stokes operators as limit transition operators between 
monodromy eigenbases}

We formulate the result from the title of the Subsection firstly 
in the case, when $k=1$, $n=2$, and then in the general case. 

{\bf Case $n=2$, $k=1$.}  Let 
$\lambda_i$, $i=1,2$, be the 
eigenvalues of the matrix $A(0)$. Without loss of generality 
we assume that $\lambda_1-\lambda_2\in \mathbb R_+$: one can achieve this 
by linear change of the time variable. 

We consider a deformation of (1.1), 
\begin{equation}
\dot z=\frac{A(t,\var)}{f(t,\var)}z,\ f(t,\var)=(t-\a_0(\var))(t-\a_1(\var)),\ 
f(t,0)\equiv t^2,\ A(t,0)=A(t), 
\tag{2.1}
\end{equation}
where $A(t,\var)$ and $f(t,\var)$ depend continuously on a 
parameter $\var\geq0$ so that $\a_0(\var)\neq\a_1(\var)$ for $\var>0$. 
Without loss of generality we assume that $\a_0+\a_1\equiv0$. 
We formulate the statement from the title of the Subsection for a generic 
deformation (2.1), see the following Definition.

\begin{definition} A family of quadratic polynomials $f(t,\var)$ 
depending continuously on a nonnegative parameter $\var$, $f(t,0)\equiv t^2$, 
with roots $\a_i(\var)$, $i=0,1$, $\a_0+\a_1\equiv0$, is said to 
be {\it generic}, if $\a_0(\var)\neq\a_1(\var)$ for $\var\neq0$, 
and the line passing through $\a_0(\var)$ and $\a_1(\var)$ 
intersects the real axis by angle bounded away from 0 uniformly in $\var$. A 
family (2.1) of {\it linear equations} with $n=2$, $k=1$, 
$\lambda_1-\lambda_2\in\mathbb R_+$ is said to be {\it generic}, 
if so is the corresponding family of polynomials $f(t,\var)$. 
\end{definition}
\begin{definition}(see, e.g., [2]). 
A singular point $t_0$ of a linear analytic ordinary 
differential equation $\dot z=\frac{B(t)}{t-t_0}z$ is said to be 
{\it Fuchsian}, if it is a first order pole of the right-hand side 
(i.e., the corresponding matrix function $B(t)$ is holomorphic at $t_0$). 
The {\it characteristic numbers} of a Fuchsian singularity 
 are the eigenvalues of the corresponding residue matrix $B(t_0)$ 
(which are equal to the logarithms divided by $2\pi i$ of the eigenvalues 
of the corresponding monodromy operator).
\end{definition}

\begin{remark} A family (2.1) of linear equations is generic, if and only if 
the difference of the characteristic numbers at $\a_0(\var)$ (or equivalently, 
at $\a_1(\var)$) of the perturbed equation is not real for small $\var$ 
and moreover has argument bounded away from $\pi\mathbb Z$ uniformly in 
$\varepsilon$ small enough. The latter condition implies that the monodromy 
operator of the perturbed equation around each singular point $\a_i$ 
has distinct eigenvalues (moreover, their modules are distinct), 
and hence, a well-defined eigenbase in the solution space (for small $\var$). 
\end{remark}

The singularities of the perturbed equation from a generic family 
have  imaginary parts of constant (and opposite) signs (by definition). 
Without loss of generality everywhere below we consider that 
$$\im\a_0>0,\ \im\a_1<0,\ \text{see Fig.2}.$$ 

\begin{figure}[ht]
  \begin{center}
    \epsfig{file=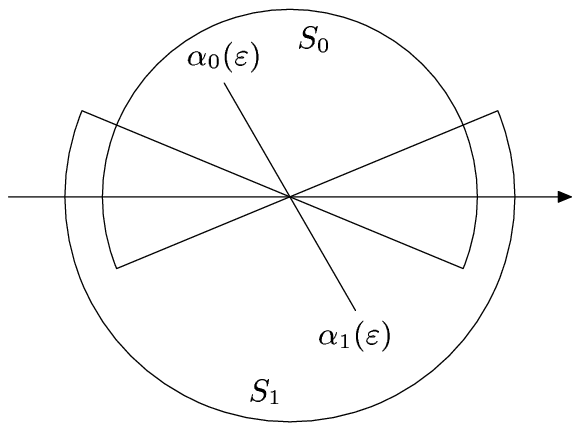}
    \caption{}
    \label{fig:2}
  \end{center}
\end{figure}

\begin{definition} Let (2.1) be a generic family of linear 
equations (see the previous Definition) whose singularity families 
satisfy the previous inequalities. 
Let $S_j$, $j=0,1$, be a pair of good sectors in the $t$- line 
(see Definition 1.2) such that for any $\var$ small enough  
$\a_j(\var)\in S_j$, $j=0,1$, $i \mathbb R_+\subset S_0$, 
$i \mathbb R_-\subset S_1$ (see Fig. 1). The sector $S_j$ is said to be 
{\it the sector associated to the singularity family} $\a_j$, $j=0,1$.  
\end{definition}
We show that appropriate 
branches of the eigenfunctions of the monodromy operator $M_i$ around $\a_i$ 
of the perturbed equation converge to canonical solutions 
of the nonperturbed equation in the corresponding sector $S_i$. This 
will imply the statement from the title of the Subsection. 

To formulate the latter statement precisely, consider the auxiliary domain  
\begin{equation}
S_i'=S_i\setminus[\a_0(\var),\a_1(\var)],\tag{2.2}\end{equation}
which is simply-connected, 
and  the canonical branches of the monodromy eigenfunctions on the domain 
$S_i'$. In more details, 
consider a small circle going around $\a_i$ and take a base point 
on it outside the segment $[\a_0(\var),\a_1(\var)]$. 
In the space of local solutions of the perturbed equation at the 
base point consider the monodromy operator $M_i$ acting by the analytic 
extension of a solution along the circle from the base point to itself 
in the counterclockwise direction. The eigenfunctions of $M_i$  
have well-defined branches (up to multiplication by constants) 
in the corresponding disc with the segment 
$[\a_0(\var),\a_1(\var)]$ deleted. Their immediate analytic extension 
yields their canonical branches on $S_i'$. In other terms, we 
identify the space of local solutions with the space of solutions on $S_i'$ 
by immediate analytic extension, consider $M_i$ as an operator acting in 
the latter space and take its eigenfunctions.

The canonical basic solutions of the nonperturbed equation 
are numerated by the indices  1 and 2, which correspond to the eigenvalues 
$\la_1,\la_2$ of $A(0)$. To state the results previously mentioned, 
let us define the analogous numeration of the monodromy eigenfunctions at 
$\a_i(\var)$. The monodromy eigenfunctions are numerated by the characteristic 
numbers (see Definition 2.2) of the corresponding singularity. The latters 
are proportional 
to the eigenvalues of the matrix $A(\a_i(\var),\var)$, which tend to $\lambda_1$ and 
$\lambda_2$, as $\var\to0$. This induces the numeration of the monodromy 
eigenfunctions by the indices 1 and 2 corresponding to the limit eigenvalues 
$\la_1$ and $\la_2$.

\begin{theorem} (see [6]).  Let (2.1) be a generic family of linear ordinary 
differential equations (see Definition 2.1), 
$\a_i(\var)$ be its singularity family, $S_i$ be the corresponding sector 
(see the previous Definition), $S_i'$ be the domain (2.2). Consider 
the eigenbase on $S_i'$ of the monodromy operator of the 
perturbed equation around $\a_i(\var)$. The appropriately normalized 
eigenbase 
(by multiplication of the basic functions by constants) converges to the 
canonical solution base (1.5) on $S_i$ of the nonperturbed equation.  
\end{theorem}

\begin{corollary} (see [6]).  Let (2.1) be a generic linear equation family 
(see Definition 2.1), $\a_i$ be its singularity families, 
$S_i$ be the corresponding sectors (see the previous Definition) chosen to 
cover a punctured neighborhood of zero, $S_i'$ be the corresponding domains 
(2.2). Let $C_0$, $C_1$ be the corresponding 
Stokes matrices (1.6) of the nonperturbed equation in the left (respectively, 
right) component of the intersection $S_0\cap S_1$. Consider 
the eigenbase on $S_i'$ of the monodromy operator of the perturbed equation 
around $\a_i(\var)$. Denote by 
$Z^i_\var(t)$  the fundamental matrix of this eigenbase. 
Let $C_0(\var)$ ($C_1(\var)$) be the transition matrix 
between the monodromy eigenbases $Z^i_\var(t)$, $i=0,1$, in the left 
(respectively, right) component of the intersection $S_0'\cap S_1'$:
\begin{equation}
Z^1_\var(t)=Z^0_\var(t)C_0(\var)\ \text{for}\ \re t<0;\ \ 
Z^0_\var(t)=Z^1_\var(t)C_1(\var)\ \text{for} \ \re t>0.\tag{2.3}
\end{equation}
For any $j=0,1$ and appropriately normalized monodromy eigenbases $Z^i_\var$, 
$i=0,1$ (the normalization of $Z_\var^0$ (only) 
depends on the choice of $j$) $C_j(\var)\to C_j$, as $\var\to0$. 
\end{corollary}

\

{\bf Case $k=1$, $n$ is arbitrary}. To state the analogues of 
 Theorem 2.5 and  Corollary 2.6  in this more general case, let us firstly extend 
the notions of 
a generic family of linear equations  and a sector associated to a singularity 
family. 

\begin{definition} Let $n,k\in\mathbb N$, $n\geq2$, 
$\Lambda=(\lambda_1,\dots,\lambda_n)$ be a set of $n$ distinct complex numbers, 
$\lambda_i\neq\lambda_j$ be a pair of them. 
A ray in $\mathbb C$ starting at 0 is called a $(k,\Lambda)$- 
{\it real dividing ray} associated to the pair $(\lambda_i,\lambda_j)$, 
if for any $t$ lying in this ray 
$\im\frac{\lambda_i-\lambda_j}{t^k}=0$. (Or equivalently, it is a ray bisecting 
an angle between two neighbor imaginary dividing rays (see Definition 1.3) 
associated to $(\lambda_i,\lambda_j)$.) 
\end{definition}

\begin{definition} Let (1.1) be an irregular equation with $k=1$, 
$\Lambda$ be the 
vector of eigenvalues of the corresponding matrix $A(0)$. Let (2.1) be its 
deformation 
depending continuously on a nonnegative parameter $\var$, $f(t,0)\equiv t^2$, 
 $\a_0+\a_1\equiv0$. The family (2.1) is said to 
be {\it generic}, if $\a_0(\var)\neq\a_1(\var)$ for $\var\neq0$, 
and the line passing through $\a_0(\var)$ and $\a_1(\var)$ 
intersects each $(k,\Lambda)$- real dividing ray by angle bounded away 
from 0 uniformly in $\var$. 
\end{definition}

\begin{definition} Let (2.1) be a generic family (see the previous Definition), 
$\Lambda$ be the corresponding eigenvalue tuple of $A(0)=A(0,0)$. Let   
$\a_0$, $\a_1$ be the corresponding singular point families, $V_i$ be the half-plane 
(depending on $\var$) 
containing $\a_i$ and bounded by the symmetry line of the segment $[\a_0,\a_1]$. 
The {\it sector associated to} $\a_i$ is a $(1,\Lambda)$- good sector 
(see Definition 1.3) independent on $\var$ that contains $V_i$ for any $\var$ small 
enough.  
\end{definition}

\begin{remark} In the previous Definition the sectors $S_0$, $S_1$ associated to 
$\a_0$, $\a_1$ respectively 
cover a punctured neighborhood of zero, so, the nonperturbed equation has a 
pair of Stokes operators ($C_0$, $C_1$) associated to this covering. 
\end{remark}

\begin{theorem} (see [6]). Let (2.1) be a generic family of linear equations (see Definition 2.8), 
$S_0$, $S_1$, $S_0'$, $S_1'$  be respectively the  corresponding associated sectors 
(see the previous Definition) and the domains (2.2). 
%corr
Then the statement of the previous Theorem remains valid.  The same is true for 
Corollary 2.6. In more details, 
consider the eigenbases $Z^i_\var$ on $S_i'$ of the monodromy operators around 
the singular points $\a_i(\var)$ of the perturbed equation, $i=0,1$. 
Let $C_0(\var)$, $C_1(\var)$ be the transition matrices (2.3) between them 
in the connected components of the intersection $S_0'\cap S_1'$. Let 
$C_0$, $C_1$ be the Stokes matrices of the 
nonperturbed equation in the corresponding limit 
connected components of the intersection $S_0\cap S_1$. Then for any 
$j=0,1$ and appropriately normalized monodromy eigenbases $Z^i_\var$ 
(the normalization of $Z^0_\var$ (only) depends on the choice of $j$) \ 
$C_j(\var)\to C_j$,  as $\var\to0$. 
%corr
\end{theorem}

\subsection{Stokes operators as limits of commutators of appropriate 
powers of the monodromy operators}

The Stokes and monodromy operators act in different linear spaces: 
in the solution spaces of the nonperturbed (respectively, perturbed) equations. 
 To formulate the statement from the title of the Subsection, let 
us firstly identify these solution spaces and specify the loops defining the 
monodromy operators. 

\def\hto{H_{t_0}}

Let (2.1) be a generic family of linear equations (in the sense of some of Definitions 
2.1 and 2.8). Take a "base point" $t_0$ in the unit disc punctured at 0.  

\begin{remark} The space of local solutions of a linear equation at a 
nonsingular point $t_0\in\mathbb C$ is identified 
with the space of initial conditions at $t_0$ (which is common for 
the nonperturbed and the perturbed equations). This identifies 
the solution spaces of the latters. The space thus obtained will be denoted by 
$H_{t_0}$.  
\end{remark}

\begin{remark} Let (1.1) be an irregular equation with 
$k=1$, $\Lambda$ be the eigenvalue tuple of the corresponding 
matrix $A(0)$. Let $S_0$, $S_1$ be $(1,\Lambda)$- good sectors 
covering a punctured neighborhood of zero in the $t$- line. 
Let $C_0$, $C_1$ be the Stokes operators (1.6) 
corresponding to the connected components of their intersection. 
Each operator $C_i$ is well-defined in the space $\hto$ 
of local solutions of (1.1) at any point $t_0$ lying in the corresponding  
component of the intersection $S_0\cap S_1$.  
\end{remark}
Now let us define the monodromy operators acting in the previous space $\hto$.  

\begin{definition} Let (2.1) be a generic family of linear equations 
(in the sense of one of Definitions 2.1 and 2.8), 
$\a_i(\var)$, $i=0,1$, be its singularity families. Fix a point 
$t_0$ (independent on $\var$) disjoint from the line passing through $\a_0(\var)$ 
and $\a_1(\var)$ for any $\var$.  
Let $l_i$ be a small circle centered at $\a_i(\var)$ whose closed disc 
is disjoint from $-\a_i(\var)$, $a_i=[t_0,\a_i]\cap l_i$, the segment 
$[t_0,a_i]$ be oriented from $t_0$ to $a_i$. Consider the closed 
path $\psi_i=[t_0,a_i]\circ l_i\circ [t_0,a_i]^{-1}$, $i=0,1$, 
which starts and ends at $t_0$ (in the case, when 
$k=1$, $n=2$, $\la_1-\la_2\in\mathbb R$, we choose $t_0\in\mathbb R$, 
see Fig.3). Define $M_i:\hto\to\hto$ to be the 
corresponding monodromy operator of the perturbed equation. 
\end{definition}

\begin{figure}[ht]
  \begin{center}
    \epsfig{file=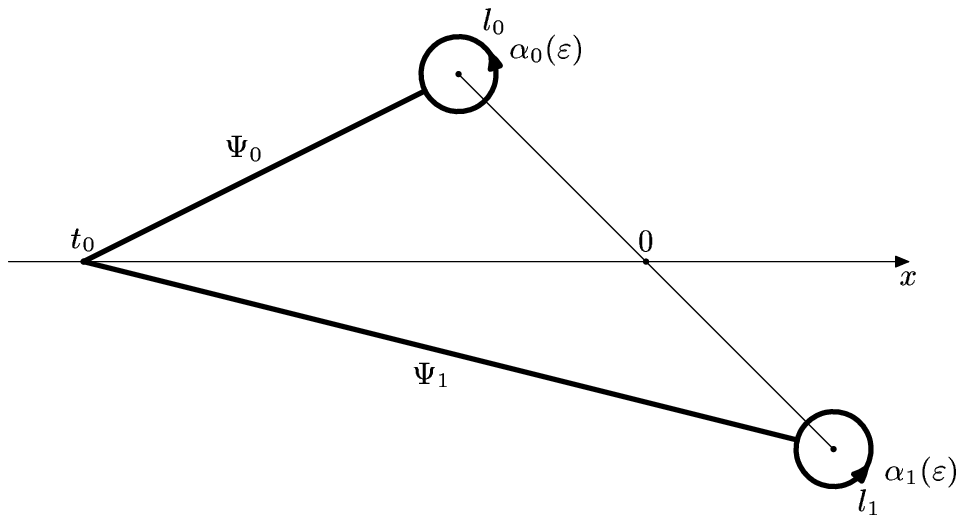}
    \caption{}
    \label{fig:3}
  \end{center}
\end{figure}

We show that commutators of appropriate noninteger powers of the operators 
$M_i$ (see the following Definition) tend to the Stokes operators. 

\begin{definition} Let $d\in\mathbb R$, $M:H\to H$ be a linear operator in a 
finite-dimensional 
linear space having distinct eigenvalues. The {\it $d$-th power of} $M$ is 
the operator having the same eigenlines, as $M$, whose corresponding 
eigenvalues are some values of the $d$-th powers of those of $M$.  
\end{definition}

Let $S_0$, $S_1$ be sectors in $\mathbb C$ with vertex at 0 covering  a 
punctured neighborhood of 0. Their {\it left (right) intersection component} is 
the component of their intersection crossed while 
going from $S_0$ to $S_1$ in the counterclockwise 
(respectively, clockwise) direction. 

\begin{theorem} Let (2.1) be a generic family of linear equations 
(in the sense of one of Definitions 2.1 and 2.8), 
$\a_i(\var)$, $i=0,1$, be its singularity families. Let $S_i$, $i=0,1$, be 
the corresponding associated sectors (see Definitions 2.4 or 2.9 respectively) 
forming a covering 
of a punctured neighborhood of zero, $C_0$, $C_1$ be the Stokes operators (1.6) 
of the nonperturbed equation corresponding to the left (respectively, right) 
component of the intersection $S_0\cap S_1$ 
%corr
(see the previous paragraph). 
%corr
Let $t_0$ be a fixed point of unit 
disc lying in the left component of the intersection $S_0\cap S_1$,  
$\hto$ be the corresponding local solution space (see Remark 2.12). 
(Then the operator $C_0$ ($C_1$) acts in the space $\hto$ (recpectively, 
$H_{-t_0}$, see 
%corr
Remark 2.13.)) 
%corr
Let 
$M_i:H_{\pm t_0}\to H_{\pm t_0}$ be the corresponding monodromy operators 
from Definition 2.14.   Then for any pair of numbers $d_0,d_1>0$ such that $d_0+d_1<1$ 
$$
M_1^{-d_1}M_0^{d_0}M_1^{d_1}M_0^{-d_0}\to C_0 \ \text{in the space}\ \hto,$$
$$M_0^{-d_0}M_1^{d_1}M_0^{d_0}M_1^{-d_1}\to C_1 \ \text{in the space}\ 
H_{-t_0},\ \text{as}\ \var\to0.$$
\end{theorem}

Theorem 2.16 is proved in  Section 3.

\def\cc{\mathbb C}
\def\oc{\overline\cc}

\subsection{The case of higher Poincar\'e rank}

Theorem 2.5 on convergence of 
the monodromy eigenbases to canonical solution bases is 
stated and proved in [6] for arbitrary 
irregular equation (for arbitrary Poincar\'e rank 
and dimension). It holds for any generic family (1.2) defined to satisfy the 
following conditions: 1) $\sum\a_i\equiv0$, $f'_\var(0,0)\neq0$ (then 
$\a_i(\var)=a_i\var^{\frac1{k+1}}(1+o(1))$, where the points $a_i$ form 
a regular polygon centered at 0); 2) no one of the previous points $a_i$ lies in a real 
dividing ray (see Definition 2.7), in other terms, no radial ray of $\a_i$ 
tends to a real dividing ray. To each singularity family $\a$ we put into 
correspondence a $(k,\Lambda)$- good sector $S$ (similarly to Definition 2.9) 
so that the canonical 
branches in $S'_\var=S\setminus\cup_{i=0}^k[0,\a_i(\var)]$ of 
the corresponding monodromy eigenfunctions converge to the 
canonical solutions of the nonperturbed equation on $S$. In the case of 
higher Poincar\'e rank for some pairs of neighbor singularities of the perturbed 
equation the corresponding sectors cannot not be chosen intersected; then 
the corresponding transition operator between the monodromy eigenbases tends  
to a product of Stokes operators. Each Stokes matrix is contained in some of the 
previous limit products, and its elements 
are expressed as polynomials in the elements of the corresponding limit product. 
On the other hand, 
%corr
in dimension two there are always two pairs of neighbor 
singularity families  such that for each singularity pair 
%corr
the corresponding sectors may be chosen intersected. Then 
the transition operator between the corresponding appropriately normalized 
monodromy eigenbases tends 
to the Stokes operator of the nonperturbed equation corresponding to the 
intersection of the sectors.

\begin{example} Consider the case, when $k=n=2$. 
%corr
Then the perturbed equation 
has three singularities, and the number of $(2,\Lambda)$- good sectors 
covering a punctured neighborhood of zero is equal to 4. 
%corr
One can prove 
the following version of Theorem 2.16. 

{\it Consider a generic deformation (1.2) of an irregular equation (1.1) 
with $k=n=2$. Let $\a_0$, $\a_1$ be a pair of singularity families 
%corr
numerated 
counterclockwise and 
%corr 
corresponding to intersected sectors (denote the latters by $S_0$ and 
$S_1$ respectively). Let $t_0\in\mathbb C\setminus0$ 
be a fixed (base) point lying between the radial rays of $\a_0(\var)$, $\a_1(\var)$ 
for all $\var$. Let $M_0$, $M_1$ be the corresponding 
monodromy operators (see Definition 2.14). Let $C$ be the 
Stokes operator corresponding to the intersection $S_0\cap S_1$. 
Then for appropriate  
$d_0,d_1\in\mathbb R\setminus0$ (depending on the family of equations)  
$$M_1^{-d_1}M_0^{d_0}M_1^{d_1}M_0^{-d_0}\to C,\ \text{as}\ \var\to0.$$
More precisely, there exist $s_i\in\mathbb N$, 
$l_0,l_1>0$ (depending on the family of equations
%corr
but not 
on $\var$) 
%corr 
such that the previous statement holds whenever $d_0$, $d_1$ 
satisfy the following system of inequalities}:
$$\begin{cases} (-1)^{s_i}d_i>0,\ i=0,1\\
 l_0d_0+l_1d_1<1\end{cases}.$$
\end{example}
\begin{remark}
The previous coefficients $l_i$  
depend on how 
close the radial rays of $\a_i$, $i=0,1$, approach the 
real dividing rays: if the minimal angle between the radial ray of $\a_i$ 
and some real dividing ray is small, then the corresponding coefficient $l_i$ 
should be chosen large enough (hence, the corresponding exponent $d_i$ 
should be taken small enough).  
\end{remark}
The author believes that the previous statements extend to the general 
case of arbitrary Poincar\'e rank and dimension. 

\section{Convergence of the commutators to Stokes operators. 
Proof of Theorem 2.16}

Firstly we prove Theorem 2.16 in the case, when $k=1$, $n=2$. Its proof 
for the case of $k=1$ and arbitrary $n$ is analogous: the modifications needed 
will be discussed in Subsection 3.4

Thus, from now on we consider that $k=1$, $n=2$, until the contrary will be 
specified. Without loss of generality we assume that 
$\Lambda=(\lambda_1,\lambda_2)=(1,-1)$. 

\subsection{Properties of the monodromy and the transition operators. 
The plan of the proof of Theorem 2.16}

Let us prove the convergence of the first commutator  from Theorem 2.16; 
the proof of the convergence of the second commutator is analogous. 

Thus, from now on we assume that the base point $t_0$ lies in the left component of the 
intersection $S_0\cap S_1$, and one can put $t_0=-\frac12$. 

\begin{definition} Consider a linear diagonalizable operator 
acting on $\mathbb C^2$ with eigenvalues of distinct modules. 
Its {\it projective multiplier} is the ratio of its eigenvalue with the 
lower module over that with the higher module. Its {\it projectivization} 
is the M\"obius transformation $\overline{\mathbb C}\to\overline{\mathbb C}$ 
induced by its action and the tautological projection $\mathbb C^2\setminus0\to
\mathbb P^1=\overline{\mathbb C}$. 
\end{definition}

\begin{remark} In the conditions of the previous Definition 
the projectivization is a hyperbolic transformation (see [1] and Definition 
4.6 in Section 4); in particular, it has an attracting fixed point. 
The projective multiplier is well-defined and its module 
is always less than 1. It is equal to the multiplier of the projectivization 
at its attracting fixed point.
\end{remark}

Let us write down the monodromy operators in the eigenbase of $M_0$ 
(which converges to the canonical solution base of the nonperturbed equation 
on $S_0$). Then the matrix of $M_0$ is diagonal: denote it 
$$\Lambda_0(\var)=\diag(\lambda_{01},\lambda_{02})(\var).$$
By Corollary 2.6, the matrix of $M_1$ is 
\begin{equation}
M_1=C(\var)\Lambda_1(\var)C^{-1}(\var), \ C(\var)\to C_0, \ \text{as}\ \var\to0,
\ \Lambda_1(\var)=\diag(\lambda_{11},\lambda_{12})(\var).\tag{3.1}
\end{equation}
The transition matrix $C(\var)$ 
tends to the Stokes matrix $C_0$, which is lower-triangular. Thus, the 
upper-triangular element of $C(\var)$ (denoted by $u(\var)$) tends to 0. 

First of all we find the asymptotics of the eigenvalues $\lambda_{ij}$ of $M_i$: 

\begin{proposition} Let (2.1) be a generic family of linear equations (see 
Definition 2.1), $t_0=-\frac12$,  $M_i$ 
be the monodromy operators of the perturbed equation from Definition 2.14, 
$f_{i1,\var}$, $f_{i2,\var}$ be their basic eigenfunctions,  
$\la_{i1}$, $\la_{i2}$ be the corresponding eigenvalues. Then 
$$\la_{01}, \la_{12}\to\infty,\ \la_{02}, \la_{11}\to0,$$
$$\ln\la_{01}=-(1+o(1))\ln\la_{02}=-(1+o(1))\ln\la_{11}=(1+o(1))\ln\la_{12},\ 
\text{as}\ \var\to0.$$
\end{proposition}
\begin{corollary} In the conditions of the previous Proposition the projective 
multipliers of $M_0$ and $M_1$ are equal respectively to 
$$\mu_0=\frac{\la_{02}}{\la_{01}},\ \mu_1=\frac{\la_{11}}{\la_{12}};\ \mu_i
\to0,\ \text{as}\ \var\to0,$$
$$ \ln\mu_0=(1+o(1))\ln\mu_1, \ \text{as}\ \var\to0.$$
\end{corollary}

\begin{proof} {\bf of Proposition 3.3.}
It follows from definition that $\ln\la_{01}=(1+o(1))\frac{2\pi i}{\a_0
-\a_1}$. The real part of the right-hand side of the previous formula is positive 
and tends to infinity (since $\im(\a_0-\a_1)>0$ by assumption, and $\a_i\to0$), 
which implies that $\la_{01}\to\infty$. The similar formulas 
written for all the $\la_{ij}$ prove the rest of the statements of the Proposition. 
\end{proof}

\def\tm{\widetilde M}
\def\tl{\widetilde{\Lambda}}

Let $d_0,d_1>0$, $d_0+d_1<1$, 
$$\tm_i=M_i^{d_i},\ \tl_i=\Lambda_i^{d_i}, \ i=0,1.$$
We prove that 
$$\tm_1^{-1}\tm_0\tm_1\tm_0^{-1}\to C_0.$$
By definition and (3.1), the matrix  of the 
previous commutator in the eigenbase of $M_0$ is 
\begin{equation} \tm_1^{-1}\tm_0\tm_1\tm_0^{-1}=C(\var)\tl_1^{-1}C^{-1}(\var)
\tl_0C(\var)\tl_1C^{-1}(\var)\tl_0^{-1}.\tag{3.2}
\end{equation}

Let $u(\var)$ be the upper-triangular element of the transition matrix 
$C(\var)$, $\mu_1(\var)$ be the projective multiplier of $M_1$. 
For the proof of the convergence to $C_0=\lim C(\var)$ of the previous commutator 
 we firstly prove that 
\begin{equation}u=O(\mu_1),\ \text{as}\ \var\to0.\tag{3.3}
\end{equation}
More precisely, we show in the next Subsection that $u=-(c_1+o(1))\mu_1$, 
where $c_1$ is the  upper-triangular element of the other Stokes matrix $C_1$. 

Let $\nu_0=\mu_0^{d_0}$, $\nu_1=\mu_1^{d_1}$ be the projective multipliers of 
the operators $\tm_0$, $\tm_1$ respectively. Formula (3.3) 
together with the previous Corollary and the condition $d_0+d_1<1$ imply 
that 
\begin{equation} u(\var)=o(\nu_0\nu_1), \ \text{as}\ \var\to0.\tag{3.4}\end{equation}

Using (3.4)  we show (in the next Lemma proved in 3.3) that if we eliminate 
subsequently the terms $C^{\pm1}(\var)$ in (3.2) (from the right to the left) 
except for the left $C(\var)$, then on each step 
the asymptotics of the modified 
expression (3.2) remains the same: the modified expression can 
%corr
be 
%corr
obtained 
from the initial one by composing it with an operator tending to the identity.  
At the last step the final modified expression will be just $C(\var)$,  
which tends to $C_0$. This will prove the convergence of (3.2) to $C_0$. 

Now the convergence of the commutator (3.2) is implied by  the 
following Lemma (modulo (3.3)).

\begin{lemma} Let $\tm_0$, $\tm_1$ be two families of two-dimensional 
%corr
complex
%corr
diagonalizable linear operators depending on a positive parameter $\var$. 
Let $\tl_i=(\la_{i1},\la_{i2})$ be the (diagonal) matrices 
of $\tm_i$, $i=0,1$, in their eigenbases. 
Let $C(\var)$ be the transition 
matrix between the eigenbases, more precisely, in the eigenbase of $\tm_0$ one 
have $\tm_1=C(\var)\tl_1 C^{-1}(\var)$. Let the eigenbases converge to some 
bases in the space so that the corresponding transition matrix $C(\var)$ tends 
to a unipotent lower-triangular matrix (denoted by $C_0$). Let 
\begin{equation} 
\nu_0=\frac{\la_{02}}{\la_{01}}\to0, \ \nu_1=\frac{\la_{11}}{\la_{12}}\to0,\ 
\text{as}\ \var\to0.\tag{3.5}\end{equation}
Let $\nu_0$, $\nu_1$ and the upper-triangular element $u(\var)$ 
of the matrix $C(\var)$ satisfy (3.4). Then 
$$\tm_1^{-1}\tm_0\tm_1\tm_0^{-1}\to C_0,\ \text{as}\ \var\to0.$$
\end{lemma}

The Lemma will be proved in Subsection 3.3.

\begin{proof} {\bf of Theorem 2.16 modulo (3.3) and Lemma 3.5.} 
The operators $\tm_i=M_i^{d_i}$ satisfy the conditions of the previous Lemma: 
all the conditions  
follow from the previous Proposition, Corollary and (3.3). This together 
with the Lemma implies the convergence of the first commutator in Theorem 2.16.  
This proves Theorem 2.16 modulo (3.3) and Lemma 3.5. 
\end{proof}

\subsection{The 
upper-triangular element of the transition matrix. Proof of (3.3)}
We prove the following more precise version of (3.3). 
\begin{lemma} Let (2.1) be a generic family of linear equations 
(see Definition 2.1), $\a_i$ be its singularity families, 
$S_i$ be the corresponding sectors (see Definition 2.4) chosen to 
cover a punctured neighborhood of zero, $S_i'$ be the corresponding domains 
from (2.2). Let $C_0$, $C_1$ be the 
Stokes matrices (1.6) of the nonperturbed equation (corresponding to the 
left (respectively, right) component of the intersection $S_0\cap S_1$), 
$$
C_0=\left(\begin{matrix} & 1 & 0 \\ & c_0 & 1\end{matrix}\right),\ 
C_1=\left(\begin{matrix} &1 & c_1\\ &0 &1\end{matrix}\right), \ \text{see Example 1.7}.
$$
Let $M_i$ be the monodromy operator of the perturbed equation around 
$\a_i(\var)$ acting in the space of solutions on $S_i'$. 
Let $Z^i_\var$  be the (fundamental matrix of) its eigenbase. 
Let $C_0(\var)$ be the transition matrix (2.3) between the bases 
$Z^i_{\var}$ in the left component of the intersection $S_0'\cap S_1'$. Let 
the previous eigenbases be normalized to converge so that $C(\var)\to C_0$ 
(see Corollary 2.6): 
$$C_0(\var)=\left(\begin{matrix} & 1+o(1) & u(\var) \\ & c_0+o(1) & 1+o(1)
\end{matrix}\right),\ u(\var)\to0.$$
Let $\la_{11},\la_{12}$ be the eigenvalues of $M_1$, 
$\mu_1=\frac{\la_{11}}{\la_{12}}$ be the corresponding projective 
multiplier. 
Then the upper-triangular element $u(\var)$ of the matrix 
$C_0(\var)$ has the asymptotics 
\begin{equation}
u(\var)=(-c_1+o(1))\mu_1^{-1},\ \text{as} 
\ \var\to0,\tag{3.6}
\end{equation}
where $c_1$ is the upper-triangular element of the Stokes matrix $C_1$. 
\end{lemma}

\begin{proof}
The transition matrix $C_0(\var)$, $Z^1_\var=Z^0_\var C_0(\var)$, compares 
the monodromy eigenbases in the left component 
of the intersection $S_0'\cap S_1'$, in particular, on a real interval 
in $\mathbb R_-$. It is not changed, when we extend the basic solutions 
analytically from $\mathbb R_-$ 
to $ \mathbb R_+$ along the real line. Denote $Z^i_{\var,+}$ the 
corresponding branch on $ \mathbb R_+$ of the extended fundamental matrix 
$Z^i_\var$, $i=0,1$. It follows from definition that 
$Z^1_{\var,+}$ is obtained from $Z^1_\var|_{\mathbb R_+}$ by applying the 
inverse monodromy operator $M_1^{-1}$: 
\begin{equation}
Z^1_{\var,+}=Z^1_\var|_{S_1'}M_1^{-1};\ \text{the matrix}\ M_1\ 
\text{is diagonal}.\tag{3.7}\end{equation}
On the other hand, we can choose a 
renormalization of the eigenbase $Z^0_{\var,+}$ by multiplication 
of the basic solutions by constants 
(i.e., changing it to $Z^0_{\var,+}L(\var)$, 
$L(\var)=diag(l_1(\var),l_2(\var))$ is some family of diagonal matrices) 
so that in the right component of the intersection $S_0'\cap S_1'$ 
the transition matrix $C_1(\var)$ between $Z^0_{\var,+}L(\var)$ 
and $Z^1_\var$ tends to the Stokes matrix $C_1$:  
$$ 
Z^0_{\var,+}L(\var)=Z^1_\var|_{S_1'}C_1(\var),\ C_1(\var)\to C_1.
$$
By definition, $Z^1_{\var,+}=Z^0_{\var,+}C_0(\var)$. Substituting the latter 
and (3.7) to the previous formula yields 
\begin{equation}
C_0(\var)=L(\var)C_1^{-1}(\var)M_1^{-1}.\tag{3.8}
\end{equation}
The matrices $C_i(\var)$ tend to the Stokes matrices $C_i$, which are unipotent. 
The matrices $L(\var)$, $M_1$ are diagonal and depend on $\var$. 
This implies that 
$$L(\var)=M_1(1+o(1)),\ \text{as}\ \var\to0.$$
This together with (3.8) implies (3.6).
\end{proof}
\subsection{Commutators\ \ of\ \ operators\ \ with\ \ asymptotically\ common\ eigenline. 
\ \ Proof of Lemma 3.5}
In the proof of Lemma 3.5 we use the following Proposition. 

\begin{proposition} Let $C(\var)$ be a family of two-dimensional 
matrices depending on a parameter $\var\geq0$ and converging to 
a unipotent lower-triangular matrix, as $\var\to0$. Let $u=u(\var)$ be the 
upper-triangular element of $C(\var)$ (thus, $u(\var)\to0$). Let 
$\Lambda(\var)=\diag
(\la_1(\var),\la_2(\var))$ be a family of diagonal matrices depending on 
$\var>0$ such that 
\begin{equation} \nu=\frac{\la_1}{\la_2}\to0, \ u=o(\nu),\ \text{as}\ \var\to0.
\tag{3.9}\end{equation}
Then 
\begin{equation}\Lambda^{-1}(\var)C(\var)\Lambda(\var)\to Id,\ \text{as}\ \var\to0.
\tag{3.10}\end{equation}
\end{proposition}
\begin{proof} The diagonal elements of the matrix in (3.10) are equal to those of 
$C(\var)$, and thus, tend to 1. Its lower-triangular element tends to 0: it 
is equal to that of $C(\var)$ (which tends to a finite limit) times $\nu$ 
(which tends to 0 by (3.9)). Its upper-triangular element, which is equal to 
$u\nu^{-1}$, tends to 0 by (3.9). This proves (3.10). 
\end{proof}

Consider the commutator (3.2): 
\begin{equation}C(\var)\tl_1^{-1}C^{-1}(\var)
\tl_0C(\var)\tl_1C^{-1}(\var)\tl_0^{-1}.\tag{3.11}\end{equation}
Using the previous Proposition, we firstly 
"kill" the right $C^{-1}(\var)$: we show that expression (3.11) is equal to 
\begin{equation}C(\var)\tl_1^{-1}C^{-1}(\var)
\tl_0C(\var)\tl_1\tl_0^{-1}(Id+o(1)).\tag{3.12}\end{equation} 
Then we kill similarly the right $C(\var)$ in (3.12) 
and the remaining $C^{-1}(\var)$.  
Finally we get that the initial commutator is equal  to $C(\var)$ times the commutator 
of diagonal matrices (which is identity) times $(Id+o(1))$. This implies that 
(3.2) tends to $C_0=\lim C(\var)$. 

The first step: killing of $C^{-1}$.  Let 
$$Q(\var)=\tl_0C^{-1}\tl_0^{-1}.$$ 
By definition, expression (3.11) is equal to 
$$C(\var)\tl_1^{-1}C^{-1}(\var)
\tl_0C(\var)\tl_1\tl_0^{-1}Q(\var).$$ It suffices to 
show that $Q(\var)\to Id$. This follows from the previous Proposition applied to 
the families of matrices $C^{-1}(\var)$ and $\Lambda(\var)=\tl_0^{-1}(\var)$:  
these families 
satisfy the conditions of the previous Proposition. Indeed, by (3.5), 
$\nu=\nu_0\to0$. The upper-triangular 
element of $C^{-1}$ (denoted by $\tilde u$) is $\tilde u=O(u)=
o(\nu_0\nu_1)=o(\nu_0)$ by (3.4). This proves 
(3.9) for $\tilde u$. The conditions of the Proposition are checked. 
Thus, by (3.10), $Q(\var)\to Id.$

The second step: killing of the right $C$ in (3.12). It repeats the previous 
 discussions with the families $C(\var)$ and $\Lambda(\var)=\tl_1\tl_0^{-1}$. 

The third step: killing of the left $C^{-1}$. Done analogously by applying 
the previous Proposition to  the matrix families $C^{-1}$ and $\tl_1$.  
Lemma 3.5 is proved. 

\subsection{Convergence of the commutators to Stokes operators: 
the higher-dimensional case}

The proof of Theorem 2.16 in higher dimensions repeats that in the two-dimensional 
case with some changes specified below. 

Let (2.1) be a generic family of equations, $\a_i$ be its singularity 
families, $i=0,1$, $S_i$, be the corresponding associated 
sectors. Consider their "left intersection component" that is crossed while going 
counterclockwise from $S_0$ to $S_1$. 
Let  $t_0$ be a point lying in this component. Let 
$\hto$ be the corresponding local solution space,   
$M_i:\hto\to \hto$  be the corresponding monodromy operators (see 
Definition 2.14). 
Let $Z^0_{\var}$ be the eigenbase of the monodromy operator $M_0$, where the 
eigenfunctions are taken in the order of decreasing of the modules of the 
corresponding eigenvalues (it appears that these modules are really 
distinct, see the next Proposition). 
Let $Z^1_{\var}$ be that of $M_1$, and the order of the eigenfunctions coincide with 
the order of increasing of the modules of the eigenvalues. 
Let $C(\var)$ be the transition matrix between them: 
$$Z^1_{\var}=Z^0_{\var}C(\var).$$
Let $C_0$, $C_1$ be the Stokes matrices of the nonperturbed equation 
in the left (respectively, right) connected component of the intersection 
$S_0\cap S_1$. 

\begin{proposition} In the above conditions for any $\var$ small enough 
 each monodromy operator $M_i$, 
$i=0,1$, has distinct eigenvalues (denote them $\la_{i1},\dots,\la_{in}$). 
Moreover, for any $j,k=1,\dots,n$, $j<k$, one has 
$\frac{\la_{0j}}{\la_{0k}}\to\infty$, $\frac{\la_{1j}}{\la_{1k}}\to0$, 
as $\var\to0$. 
Appropriately normalized eigenbases $Z^0_{\var}$ and $Z^1_{\var}$ converge 
to canonical solution bases of the nonperturbed equation 
in $S_0$ and $S_1$ respectively.  Let 
$\la_1,\dots\la_n$ be the eigenvalues of the matrix $A(0)$ numerated in 
the order of increasing of the values $\re\frac{\la_j}{i\a_0(\var)}$. 
The numeration of each (converging) monodromy eigenbase corresponds to the 
numeration of the limit canonical solution base by the previous eigenvalues 
$\la_j$. 
\end{proposition}

In dimension two the Proposition follows from Proposition 3.3. 
In higher dimension its proof is analogous to that of Proposition 3.3. 

The Stokes matrices $C_0$ and $C_1$ are lower- (respectively, upper-) 
triangular. This is implied by the last statement of the previous Proposition 
and the following well-known fact. 
\begin{proposition} (see, e.g., [9]). Let $k,n\in\Bbb N$, $n\geq2$, \  
(1.1) \ be\ an irregular equation, 
$\Lambda=(\la_1,\dots,\la_n)$ be the eigenvalues of the corresponding 
matrix $A(0)$. Let $S_0$, $S_1$ be a pair of intersected $(k,\Lambda)$- good 
sectors. Let there exist a $t\in S_0\cap S_1$ such that the sequence of 
the values $\re\frac{\la_j}{t^k}$, $j=1,\dots,n$, increases. 
Consider the canonical sectorial solution bases 
of (1.1) numerated by $\la_j$. 
Then the Stokes matrix of (1.1) in the connected 
component containing $t$ of the intersection $S_0\cap S_1$ is lower-triangular. 
In the case of the inverse order of the eigenvalues it is upper-triangular.
\end{proposition}

The point $t=i\a_0(\var)$ satisfies the conditions of the previous Proposition 
with $k=1$ (the last statement of Proposition 3.8). Hence, by Proposition 3.9, 
the Stokes matrix $C_0$ is lower-triangular and $C_1$ is upper-triangular. 

Let 
$$ C(\var)= (C_{ij}(\var)),\ C_1'=C_1^{-1}=(C'_{1,ij}).$$
Formula (3.6) of Lemma 3.6 extends to higher dimension as follows: 
\begin{equation}C_{jk}(\var)=(C'_{1,jk}+o(1))(\frac{\la_{1j}}{\la_{1k}}),\ 
\text{as}\ \var\to0.\tag{3.13}
\end{equation} 
The proof of (3.13) repeats that of (3.6) in Subsection 3.2. 

Let $d_0,d_1>0$, $d_0+d_1<1$, 
$\widetilde\la_{ij}=\la_{ij}^{d_i}$, $i=0,1$, $j=1,\dots,n$. One has 
\begin{equation} C_{jk}(\var)=
O(\frac{\widetilde\la_{1j}}{\widetilde\la_{1k}}
\frac{\widetilde\la_{0k}}{\widetilde\la_{0j}}),\ \text{as}\ \var\to0.
\tag{3.14}
\end{equation} 
Formula (3.14) follows from (3.13),  the inequality $d_0+d_1<1$ 
and the asymptotic formula 
$\ln\la_{0j}=-(1+o(1))\ln\la_{1j}$, $j=1,\dots,n$, which is proved 
analogously to Proposition 3.3. 
As at the end of Subsection 3.1, Theorem 2.16 is implied 
by (3.14) and the following higher-dimensional analogue of Lemma 3.5. 

\begin{lemma} Let $\tm_0$, $\tm_1$ be two families of $n$- dimensional 
diagonalizable linear operators depending on a positive parameter $\var$. 
Let $\tl_i=\diag(\la_{i1},\dots,\la_{in})$  be the  (diagonal) matrices 
of $\tm_i$ in their eigenbases. Let $C(\var)=(C_{jk}(\var))$ be the transition 
matrix between their eigenbases, more precisely, in the eigenbase of $\tm_0$   
the matrix of $\tm_1$ is $C(\var)\tl_1 C^{-1}(\var)$. Let the eigenbases converge to some 
bases in the space so that the transition matrix $C(\var)$ converges to 
a unipotent lower-triangular matrix (denoted by $C_0$). Let for any 
$j<k$, $j,k=1,\dots,n$, 
$$\frac{\la_{0j}}{\la_{0k}}\to\infty, \ 
\frac{\la_{1j}}{\la_{1k}}\to0, \ 
\text{as}\ \var\to0.$$
Let in addition the asymptotic formula (3.14) hold. Then 
$$\tm_1^{-1}\tm_0\tm_1\tm_0^{-1}\to C_0,\ \text{as}\ \var\to0.$$
\end{lemma}
The proof of the Lemma repeats that of Lemma 3.5 with obvious changes. 

\section{Generic divergence of  monodromy operators along 
degenerating loops}

In the present Section we consider only 
{\it two-dimensional irregular equations with Poincar\'e rank $k=1$} 
and their generic deformations. As before, without loss of generality we 
assume that $\la_1-\la_2\in\mathbb R_+$, $\im\a_0>0$, $\im\a_1<0$.

Let (2.1) be a generic family of linear equations, $\a_i$, $i=0,1$, 
be its singularity families, $S_0$, $S_1$ be the corresponding 
associated sectors forming a covering of a punctured neighborhood of 0. Let 
$t_0\in \mathbb R_-$ be arbitrary fixed base point independent on $\var$. 
Let  $M_0=M_0(\var)$, $M_1=M_1(\var)$ be the 
corresponding monodromy operators of the perturbed equation 
(see Definition 2.14). 

Consider the circle centered at 0 and passing through $t_0$ with the 
counterclockwise orientation (it bounds a disc containing 
both singularities of the perturbed equation for any small $\var$). 
The monodromy operator along the previous circle 
is called the {\it complete monodromy}. 

\begin{remark} The complete monodromy of the perturbed equation in a generic 
family (2.1) converges to the monodromy of the nonperturbed equation along 
the counterclockwise circuit. In the previous conditions the complete monodromy 
is equal to $M_0M_1$. 
\end{remark}

In the present Section we state and  prove the Theorem saying 
that for any generic deformation (2.1) of a typical equation (1.1) (see the next 
Definition)  
 each word (1.3) with integer exponents tends to infinity in $GL_n$, except for 
the powers $(M_0M_1)^k$ of the complete monodromy. 

\subsection{The statement of the divergence Theorem}

\begin{definition} Let (1.1) be an irregular equation, as at the beginning of 
the paper, $t_0\in\cc\setminus0$ be arbitrary fixed base point, 
$M:\hto\to\hto$ be the counterclockwise monodromy operator around zero.  
Consider some branches at $t_0$ of all the sectorial 
canonical solutions of (1.1) as elements of $\hto$ 
and take the collection of the complex lines in $\hto$ generated by them. 
The equation is said to be 
{\it typical}, if for any $k\in\mathbb Z\setminus0$ 
no line from the previous collection is transformed by $M^k$ to another line 
from the same collection.  
\end{definition}

\begin{remark} The definition of typical equation does not depend on the choice of 
the base point and the branches of the canonical solutions. 
The condition that an equation (1.1) is typical is equivalent to a countable 
number of polynomial inequalities on the formal monodromy eigenvalues and the 
elements of the Stokes matrices. 
\end{remark} 

\begin{remark} Both Stokes operators of a typical equation are nontrivial. 
\end{remark}

Each monodromy operator word (1.3) can be rewritten as 
\begin{equation}
M_{j_n}^{s_n}\dots M_{j_1}^{s_1}, \ s_i=\pm 1,\ M_{j_k}^{s_k}
M_{j_{k-1}}^{s_{k-1}}\neq1\ \text{for any}\ k=2,\dots,n,\tag{4.1}
\end{equation} 
where $n=\sum_{i=1}^l|d_i|$ in the notations of (1.3). We consider those words 
(4.1) that do not coincide literally with powers of the complete monodromy.  

\def\tm{\widetilde M}
\def\tmm{\widetilde m}
\begin{definition} A word (4.1) 
is said to be {\it reduced}, if it does not coincide (literally) with 
neither $M_0M_1\dots M_0M_1$, nor $M_1^{-1}M_0^{-1}\dots M_1^{-1}M_0^{-1}$. 
\end{definition}
\begin{theorem} Let (1.1) be a typical equation (see Definition 4.2), 
(2.1) be its generic deformation, $t_0$, $M_0,M_1:\hto\to\hto$, be as at 
the beginning of the Section. Then any monodromy operator given by 
a reduced word (4.1) tends to infinity (together with its 
projectivization, see Definition 3.1), as $\var\to0$.  
\end{theorem} 
\subsection{Projectivization. The scheme of the proof of Theorem 4.6}
 Instead of  invertible linear operators $\mathbb C^2\to\mathbb C^2$ 
we will consider their projectivizations, which are M\"obius transformations 
$\overline{\mathbb C}\to\overline{\mathbb C}$. 
Denote $m_0$, $m_1$ the projectivizations of the monodromy operators 
$M_0,M_1:\mathbb C^2\to\mathbb C^2$. We show that any reduced word 
$\tmm=m_{j_n}^{s_n}\dots m_{j_1}^{s_1}$ tends to infinity in 
the M\"obius group, as $\var\to0$. This will prove the Theorem. 

Recall the following 
\begin{definition} (see [1]). A M\"obius transformation 
$m:\overline{\mathbb C}\to
\overline{\mathbb C}$ is said to be {\it hyperbolic}, if it has one repelling 
fixed point (then there is a unique attracting fixed point and each orbit except for 
the repeller tends to the attractor). A hyperbolic transformation with 
repeller $a$ and attractor $b$ will be presented as the picture at Fig.4 
\end{definition}

\begin{figure}[ht]
  \begin{center}
    \epsfig{file=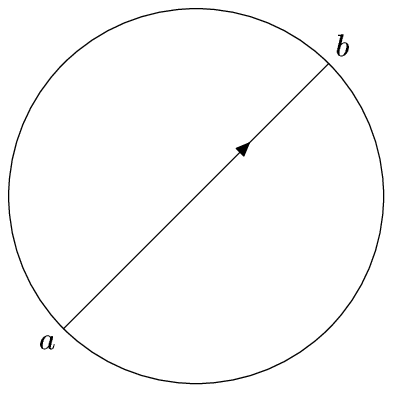}
    \caption{}
    \label{fig:4}
  \end{center}
\end{figure}

In the proof of the divergence of a reduced word $\tmm$ of the 
projectivizations we use their following properties. 

\begin{proposition} Let (2.1), $S_i$, $t_0$, 
$M_i$ be as at the beginning of the Section (the nonperturbed equation is not 
necessarily typical),  
$m_i$ be the projectivizations of $M_i$, $i=0,1$ (see Definition 3.1). 
Let $Z^i_\var=(f_{i1,\var},f_{i2,\var})$ be the 
eigenbases of $M_i$, $(f_{_i1}, f_{i2})$ be the sectorial canonical solution 
bases on $S_i$ of the nonperturbed equation, $i=0,1$. Let 
$p_{ij,\var}, p_{ij}\in\overline{\mathbb C}$ be the tautological projection 
images of $f_{ij,\var}$ and $f_{ij}$ respectively. Then $m_i$ are hyperbolic 
transformations (see the previous Definition) with fixed points 
$p_{ij,\var}$: $p_{02,\var}$, $p_{01,\var}$ are respectively the repelling and 
attracting fixed points of $m_0$; $p_{11,\var}$, $p_{12,\var}$ are respectively 
the repelling and attracting fixed points of $m_1$. 
\end{proposition}

\begin{proof} The Proposition follows from Proposition 3.3
\end{proof}
\begin{proposition} Let (2.1), $S_i$, $t_0$, $M_i$ be as at the beginning of the 
Section, $p_{ij}$, $p_{ij,\var}$ be 
the tautological projection images 
of the canonical basic solutions of the nonperturbed equation and 
the eigenfunctions of $M_i$ respectively (see the previous Proposition). Then 
\begin{equation}p_{02}=p_{12},\ p_{ij}=\lim_{\var\to0}p_{ij,\var},\ 
\text{see Fig.5a,b}.\tag{4.2}
\end{equation}
\end{proposition}

\begin{figure}[ht]
  \begin{center}
    \epsfig{file=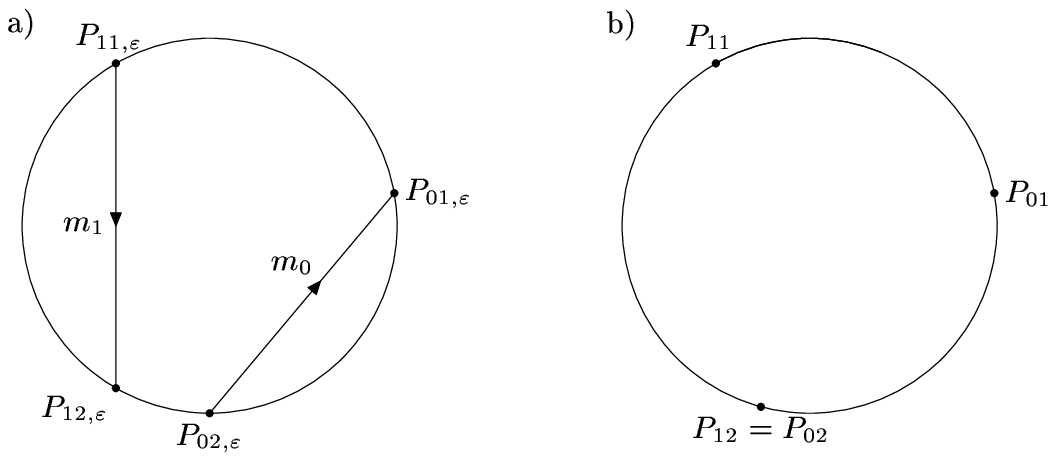}
    \caption{}
    \label{fig:5}
  \end{center}
\end{figure}

\begin{proof} The statements on the limits in (4.2) follow from Theorem 2.5.  
The coincidence of $p_{02}$ and $p_{12}$ in (4.2) follows from the 
lower-triangularity of the Stokes matrix $C_0$ (see Example 1.7).
\end{proof}

As it is shown below, Theorem 4.6 is implied by the two previous Propositions 
and the following Lemma. 

\begin{lemma} Let $p_{02}=p_{12}$, $p_{01}$, $p_{11}$ be a given triple 
of distinct points in $\overline{\mathbb C}$. Let 
$m_i=m_i(\var)$, $i=0,1$, be two families of hyperbolic M\"obius 
transformations depending on a positive parameter $\var$, 
$\mu_0$, $\mu_1$ be the multipliers of their attractors. 
Let $p_{01,\var}$, $p_{12,\var}$ be respectively 
the attractors of $m_0$ and $m_1$, $p_{02,\var}$, $p_{11,\var}$ be their 
repellers. Let $p_{ij,\var}\to p_{ij}$, $\mu_i\to0$, as $\var\to0$. 
Let in addition the product $m_0m_1$ converge to a M\"obius transformation 
$m$ such that for any $k\in\mathbb Z\setminus0$, $i=0,1,$ $j=1,2$ 
the image $m^kp_{ij}$ coincide with no other $p_{ls}$. 
Then any reduced word $\tmm=m_{j_n}^{s_n}\dots m_{j_1}^{s_1}$ 
(see Definition 4.5) tends to infinity.
\end{lemma}

The 
projectivizations of the monodromy 
operators satisfy the conditions of the Lemma. Indeed, the multipliers 
$\mu_i$ tend to 0 by Corollary 3.4. The convergence $p_{ij,\var}\to p_{ij}$ 
follows from the previous Proposition. 
The product $m_0m_1$ tends to the projectivization (denoted by $m$) 
of the monodromy of the nonperturbed equation. The inequalities $m^kp_{ij}\neq p_{ls}$ 
follow from typicality (see Definition 4.2). This together with the Lemma 
proves divergence of $\tmm$. Theorem 4.6 is proved modulo the Lemma.

\subsection{Divergence of word of projectivizations. Proof of Lemma 4.10}
 As it is shown below, Lemma 4.10 is implied by the following statement. 

\begin{lemma} In the conditions of the previous Lemma let  
$\tmm=m_{j_n}^{s_n}\dots m_{j_1}^{s_1}$ be a reduced word such that 
\begin{equation} m_{j_n}^{s_n}m_{j_{n-1}}^{s_{n-1}}\neq(m_0m_1)^{\pm1}.\tag{4.3}\end{equation}
Let $x\in\overline{\mathbb C}$ be arbitrary point such that 
\begin{equation} m^sx\neq p_{ij}\ \text{for any} \ 
s=-n,\dots,n, \ i=0,1, \ j=1,2;\ m=\lim_{\var\to0}m_0m_1.\tag{4.4}\end{equation}
Then the image $\tmm x$ converges to the limit of the 
attractor of $m_{j_n}^{s_n}$. 
\end{lemma}
Let us prove Lemma 4.10 modulo Lemma 4.11. If the reduced word $\tmm$ under consideration 
satisfies (4.3), then it tends to infinity. Indeed, by Lemma 4.11, 
the image $\tmm x$ of a generic $x$ tends to the attractor of $m_{j_n}^{s_n}$, 
hence, $\tmm\to\infty$. Otherwise, $\tmm=(m_0m_1)^km'$, where 
$k\in\mathbb Z\setminus0$, $m'$ is a word satisfying (4.3) 
of a length less than that of $\tmm$. The new word $m'$  
tends to infinity by the previous statement. The product $m_0m_1$ in 
the previous expression for $\tmm$ has a finite limit. Hence, $\tmm$ tends 
to infinity as well. Lemma 4.10 is proved.

\begin{proof} {\bf of Lemma 4.11.} In the proof of Lemma 4.11  
we use the following obvious 

\begin{proposition} Let $m'=m'(\var)$ be a family of hyperbolic 
M\"obius transformations depending on a positive parameter $\var$, $m'\to\infty$, 
as $\var\to0$, so that the attractor and the repeller of $m'$ tend to 
distinct limits (hence, the multiplier of the attractor 
tends to 0). Then for any point $x\in\overline{\mathbb C}$ distinct from 
the limit of the repeller its image $m'x$ tends to the same limit, as the 
attractor. The same statement 
holds for arbitrary family $x(\var)$ of points bounded away from the limit of 
the repeller. 
\end{proposition}

We prove Lemma 4.11 by induction in the 
length $n$ of the word. For $n=1$ its statement is obvious. Suppose we 
have proved the Lemma for the words of any length less than a given $n$. 
Let us prove it for a word $\tmm=m_{j_n}^{s_n}\dots m_{j_1}^{s_1}$ of the 
length $n$. 

Without loss of generality we assume that $m_{j_n}^{s_n}=m_0$: the contrary 
case is treated analogously. Then by (4.3) and the inequality 
$m_{j_n}^{s_n}m_{j_{n-1}}^{s_{n-1}}\neq1$ (see (4.1)), 
\begin{equation} m_{j_{n-1}}^{s_{n-1}}\neq m_1, m_0^{-1},\ \text{thus,}\ 
m_{j_{n-1}}^{s_{n-1}}=m_0 \ \text{or}\ m_1^{-1}.
\tag{4.5}\end{equation}
Consider the word 
$$m'=m_{j_{n-1}}^{s_{n-1}}\dots m_{j_1}^{s_1}=m_{j_n}^{-s_n}\tmm.$$
Suppose firstly that it satisfies (4.3). Then for any  $x$ satisfying (4.4) 
$m'x$ tends to the limit (denoted $p_{ij}$) of the 
attractor of $m_{j_{n-1}}^{s_{n-1}}$ by the induction hypothesis. The 
latter attractor can be either $p_{01,\var}$, or $p_{11,\var}$, which are 
the attractor of $m_0$ and the repeller of $m_1$ respectively. 
This follows from (4.5). Hence, the limit $p_{ij}$ is either $p_{01}$, or 
$p_{11}$; in both cases it does not coincide with the limit $p_{02}$ 
of the repeller 
of $m_0$. Therefore, the image $m_0p_{ij}$ (and hence, $m_0(m'x)=\tmm x$) 
tends to the same limit, as the attractor of $m_0$ 
(by the previous Proposition). 

Now suppose that the word $m'$ does not satisfy (4.3). Then 
$$m'=(m_0m_1)^km'',$$
where $k\in\mathbb Z\setminus0$, and $m''$ is a word satisfying (4.3) 
of a length less than that of $m'$. The induction hypothesis applied to 
$m''$ implies that for any $x$ satisfying (4.4) its image $m''x$ tends to 
the limit of the attractor of the last left element of the word $m''$, thus, 
to some $p_{ij}$. Therefore, $m'x\to m^kp_{ij}$. The image $m^kp_{ij}$ 
coincides with no $p_{sl}$ (in particular, with the limit $p_{02}$ of the 
repeller of $m_0$). This follows from the last condition of Lemma 4.10. 
Therefore, $m_0m^kp_{ij}$ tends to the limit of 
the attractor of $m_0$ (and so does $m_0m'x=\tmm x$ 
by the previous Proposition). The induction step is over. Lemma 4.11 is proved.
The proof of Lemma 4.10 is completed.\end{proof}

\bigskip

I am grateful to V.I.Arnold, Yu.S.Ilyashenko and A.A.Bolibrukh for 
stating the problems. I wish to thank them and V.Kleptsyn for their interest 
to the work and helpful remarks. The research was supported by part 
by CRDF grant  RM1-2358-MO-02 and by 
Russian Foundation for Basic Research (RFFI) grant 02-02-00482. 

%\section{References}

\textit{CNRS, Unit\'e de Math\'ematiques Pures et Appliqu\'ees, M.R.,} \\
\textit{\'Ecole Normale Sup\'erieure de Lyon} \\
\textit{46 all\'ee d'Italie, 69364 Lyon Cedex 07 France}
%\\[12pt]}
\end{document}